\documentclass{amsproc}

\usepackage{amsmath,amssymb}
\usepackage{amssymb,a4}

\usepackage{amsthm}

\def\z{\mathbb{Z}}

\def\ad{\hbox{ad}\,}
\def\ch{\hbox{\rm char}\,}

\def\a{\alpha}
\def\b{\beta}
\def\sg{\sigma}

\def\gl{\frak{gl}}
\newfont{\df}{eufm10}

\newcommand{\ml}[1]{{\mathcal #1}}

\def\sg{\sigma}

\def\fsl{\frak{sl}\,}
\def\fst{\frak{st}\,}
\def\stl{\frak{stl}\,}
\def\vp{\varphi}

\def\ot{\otimes}

\def\tr{\hbox{\rm tr}}
\def\str{\hbox{\rm str}}

\def\ker{\hbox{Ker}\,}

\def\ad{\hbox{\rm ad}\,}
\def\im{\hbox{\rm Im}\,}

\def\ot{\otimes}

\def\mg{{\bf \frak g}}

\newtheorem{theo}{Theorem}[section]
\newtheorem{defi}[theo]{Definition}
\newtheorem{lemm}[theo]{Lemma}
\newtheorem{coro}[theo]{Corollary}
\newtheorem{prop}[theo]{Proposition}

\theoremstyle{definition}
\newtheorem{remark}[theo]{Remark}

\numberwithin{equation}{section}

\begin{document}

\title[Universal Central Extensions of $\fsl(m, n, \ml A)$]
{Universal Central Extensions of the Matrix \\ Leibniz Superalgebras
$\fsl(m, n, \ml A)$}

\author[Hu]{Naihong Hu}
\address{Department of Mathematics, East China Normal University,
Shanghai 200062, PR China} \email{nhhu@euler.math.ecnu.edu.cn}
\thanks{N.H., supported in part by the NNSF (Grant 10431040), the
TRAPOYT and the FUDP from the MOE of China, the SRSTP from the
STCSM, the Shanghai Priority Academic Discipline from the SMEC}

\author[Liu]{Dong Liu}
\address{Department of Mathematics, East China Normal University,
Shanghai 200062, PR China}
\thanks{D.L., supported by the `Asia-Euro Link Programme'}
\curraddr{Department of Mathematics, Changzhou Institute of
Technology, Jiangsu 213022, PR China }
\email{liud@czu.cn}

\date{}


\begin{abstract}
The universal central extensions and their extension kernels of the
matrix Lie superalgebra $\fsl(m, n, \ml A)$, the Steinberg Lie
superalgebra $\fst(m, n, \ml A)$ in category {\bf SLeib} of Leibniz
superalgebras are determined under a weak assumption (compared with
\cite{MP}) using the first Hochschild homology and the first cyclic
homology group.
\end{abstract}

\date{}

\maketitle



\section{Introduction}

Leibniz algebra was introduced by Loday \cite{Lo1}, studied by
Cuvier \cite{C} and others. Loday-Pirashvili \cite{LP} established
the concept of universal enveloping algebras of Leibiniz algebras
and interpreted the Leibniz (co)homology $HL_*$ (resp. $HL^*$) as a
$\text{\it Tor}$-functor (resp. $Ext$-functor). The central
extensions of Leibniz algebras and Lie superalgebras have been
investigated recently (for instance, see \cite{Lo1}, \cite{Lo2},
\cite{Gao2}, \cite{Gao3}, \cite{IK}, \cite{LH1}, \cite{LH2},
\cite{MP}, etc.). Leibniz superalgebra and its cohomology were
further discussed by Dzhumadil'daev in \cite{D}. The universal
central extension of the matrix Lie superalgebras over an
associative algebra $\ml A$ in category {\bf SLie} of Lie
superalgebras was obtained in \cite{MP} under the assumption that
$m+n\ge 5$.
\begin{theo} \cite{MP}
If $m+n\ge5$, the universal central extension of $\fsl(m, n, \ml
A)$ in category {\bf SLie} is $\fst(m, n, \ml A)$ with kernel
$HC_1(\ml A)$, where $HC_1(\ml A)$ is the first cyclic homology
group of $\ml A$.
\end{theo}

Mainly motivated by \cite{MP}, \cite{LP}, \cite{AG}, \cite{Gao1},
\cite{Gao2} and \cite{KL}, we get the universal central extension of
the matrix Lie superalgebra $\fsl(m, n, \ml A)$ in category {\bf
SLeib} of Leibniz superalgebras and prove a main theorem under a
weak assumption that $m+n\ge 3$:
\begin{theo}
If $m+n\ge3$ with $\ch K\ne 2$ if $m+n=4$, $\ch K\ne 3$ if
$m+n=3$, the universal central extension of $\fsl(m, n, \ml A)$ in
category {\bf SLeib} is $\stl(m, n, \ml A)$ with kernel $HH_1(\ml
A)$, where $HH_1(\ml A)$ is the first Hochschild homology group of
$\ml A$.
\end{theo}

More precisely, we determine the universal central extension and its
kernel of $\fsl(m, n, \ml A)$ in category {\bf SLeib} in Sections 3,
4; and so does for the Steinberg Lie superalgebra $\fst(m, n, \ml
A)$ in {\bf SLeib} as well in Section 5.

\section{Leibniz superalgebras}

Throughout this paper, $K$ denotes a field with characteristic $\ne
2, 3$, $\ml A$ an associative unital $K$-algebra.

\begin{defi} \cite{D}  A Leibniz superalgebra is a
$\z_2$-graded $K$-vector space $L=L_{ 0}\oplus L_{ 1}$ with a
$K$-bilinear map $[-,-]: L\times L\to L$ satisfying $[L_{\sg},
L_{\sg'}]\subseteq L_{\sg+\sg'}$ ($\sg,\, \sg'\in\z_2$) and the
Leibniz identity $[[a, b], c]= [a, [b, c]]-(-1)^{|a||b|}[b, [a,
c]]$, for homogenous $\,a, \,b, \,c\in L$, where $|a|$ denotes the
degree of $\,a$ for a homogenous element $a\in L$.
\end{defi}

Clearly, $L_{0}$ is a Leibniz algebra; any Lie superalgebra is a
Leibniz superalgebra; any Leibniz algebra is a trivial Leibniz
superalgebra. A Leibniz superalgebra is a Lie superalgebra if and
only if $[\,a, b\,]+(-1)^{|a||b| }[b, a]=0$, for homogenous $a,\,
b\in L$.

For a Leibniz superalgebra $L$, let $L_{SLie}$ be the quotient of
$L$ by the ideal generated by elements $[x, y]+(-1)^{|x||y|}[y,
x]$, for homogenous $x, y\in L$. Clearly, $ L_{SLie}$ is a Lie
superalgebra. The canonical projection $\pi: L\to L_{SLie}$ is
universal among the maps from $L$ to Lie superalgebras. In fact,
the functor $(-)_{SLie}:$ {\bf SLeib}$\to${\bf SLie} is left
adjoint to $inc:$ {\bf SLie}$\to${\bf SLeib}. The cohomology of
Leibniz superalgebras has been defined in \cite{D}. The following
results are clear.

\begin{prop} \cite{LH3}
A Leibniz superalgebra $L$ admits a universal central extension
$\hat L$ if and only if $L$ is perfect $($i.e., $[L, L]=L$$)$.
\end{prop}

\begin{lemm} \cite{LH3}
Let $(X, f)$ and $(Y, g)$ be two central extensions of a Leibniz
superalgebra $L$. If $Y$ is perfect, then there exists only one
homomorphism $h$ from $Y$ to $X$ such that $f\circ h=g$.
\end{lemm}

 Consider the matrix Lie superalgebra
$$\gl(m, n, \ml A):=\left\{\,X \left| \ X=
\left(\begin{array}{c c} A & B \\  C & D
\end{array}\right)\,\right\}\right.,$$ where $A, B, C, D$ are
matrices in size $m\times m, m\times n, n\times m, n\times n$
respectively with coefficients in $\ml A$, and $m, n\ge0, m+n\ge
2$. Its supercommutator is defined as $[X, Y]=XY-(-1)^{\a\b}YX$
for $X\in \gl(m, n,\ml A)_\a$,  $Y\in \gl(m, n, \ml A)_\b$, $\a,
\b\in\z_2$.

By definition, the {\it special linear Lie superalgebra} with
coefficients in $\ml A$ is
$$
\fsl(m, n, \ml A):=\gl(m, n, \ml
A)^{(1)}=\{X\in\gl(m, n, \ml A)\mid \hbox{str}X=0\},
$$
where $\hbox{str}X:=\hbox{tr}A-\hbox{tr}D$ is the supertrace of $X$.
Note that if $n\ne m$, the Lie superalgebra $\fsl(m, n, \ml A)$ is
simple. Obviously, $\fsl(m, n, \ml A)$ has generators $E_{ij}(a)
\,(1\le i\ne j\le m+n, a\in \ml A)$, and subject to the relations
below:
\begin{eqnarray*}
&&[E_{ij}(a), E_{kl}(b)]=0,\  \hbox{if}\ i\ne l,\ \hbox{and}\  j\ne k;\\
&&[E_{ij}(a), E_{kl}(b)]=E_{il}(ab),\  \hbox{if}\ i\ne l,\ \hbox{and}\  j= k;\\
&&[E_{ij}(a), E_{kl}(b)]=-(-1)^{\tau_{ij}\tau_{kl}}E_{kj}(ba),\
\hbox{if}\ i= l,\ \hbox{and}\  j\ne k,
\end{eqnarray*}
where
\begin{equation*}
\tau_{ij}=\begin{cases}0, \quad\hbox{ if } \ 1\le i, j\le m, \
\hbox{ or } \ m+1\le i,j\le m+n;\cr 1, \quad\hbox{if } \ 1\le i\le
m, j>m, \ \hbox{ or } \  1\le j\le m, i>m.\cr
\end{cases}
\end{equation*}

For $\tau_{ij}$, it is clear that $\tau_{ij}=\tau_{ji}$ and
$$(-1)^{\tau_{ij}+\tau_{jk}}=(-1)^{\tau_{ik}}.\eqno(2.1)$$
Note that $\fsl(m, n, \ml A)=\fsl(m, n, \ml A)_{0}
\bigoplus\fsl(m, n, \ml A)_{1}$, $\fsl(m, n, \ml A)_{\a}=\langle
E_{ij}(a), a\in\ml A \mid \tau_{ij}=\a\rangle$, $\a\in\z_2$.

By definition, the {\it Steinberg Lie superalgebra} $\fst(m, n,
\ml A)$ is a Lie superalgebra generated by symbols $u_{ij}(a)$,
$1\le i\ne j\le m+n$, $a\in \ml A$, subject to the relations
\begin{gather*}
 u_{ij}(k_1a+k_2b)=k_1u_{ij}(a)+k_2u_{ij}(b),\
\hbox{ for } \ a, b\in\ml A, \ k_1, k_2\in K;\tag{1}\\
[u_{ij}(a), u_{kl}(b)]=0,\  \hbox{ if }\ i\ne l,\ \hbox{ and } \  j\ne k;\tag{2}\\
[u_{ij}(a), u_{kl}(b)]=u_{il}(ab),\  \hbox{ if } \ i\ne l,\ \hbox{and } \  j= k;\tag{3}\\
[u_{ij}(a), u_{kl}(b)]=-(-1)^{\tau_{ij}\tau_{kl}}u_{kj}(ba),\ \hbox{
if } \ i= l,\ \hbox{and} \  j\ne k.\tag{4}
\end{gather*}

Now we define the {\it Steinberg Leibniz superalgebra} $\stl(m, n,
\ml A)$ to be a Leibniz superalgebra with generators $v_{ij}(a)$,
$1\le i\ne j\le m+n$, $a\in \ml A$ and subject to the above
relations (1)---(4).

Note that $\stl(m, n, \ml A)=\stl(m, n, \ml A)_{0}\bigoplus\stl(m,
n, \ml A)_{1}$, $\stl(m, n, \ml A)_{\a}=\langle v_{ij}(a)\mid
a\in\ml A,  \tau_{ij}=\a\rangle$ for $\a\in\z_2$. Clearly,
relations (3)---(4) make sense only if $m+n\ge3$.

Define homomorphisms $\vp$ (resp. $\,\psi$) of Lie (resp. Leibniz)
superalgebras
$$\vp \ (\hbox{\it resp.} \ \psi): \quad \fst(m, n, \ml A) \quad
(\,\hbox{\it resp.} \ \stl(m, n, \ml A)\,)\longrightarrow \fsl(m, n,
\ml A)$$ as $\vp(u_{ij}(a))=E_{ij}(a)$ (\,resp.
$\psi(v_{ij}(a))=E_{ij}(a)$\,). Clearly, $\vp$, $\psi$ are
surjective. Moreover, $\fst(m, n, \ml A)=\stl(m, n, \ml A)_{SLie}$.
Denote this projection by $\pi:\stl(m, n, \ml A)\longrightarrow
\fst(m, n, \ml A)$.

\begin{lemm}
The Steinberg Leibniz superalgebra $\stl(m, n, \ml A)$ with
$m+n\ge3$ is perfect.
\end{lemm}

\section{Universal central extension of $\fsl(m, n, \ml A)$ in {\bf SLeib}}

\begin{theo}
If $m+n\ge3$, then $(\stl(m, n, \ml A), \psi)$ is a central extension of the Leibniz superalgebra $\fsl(m, n, \ml A)$.
\end{theo}

To begin with the proof, let us calculate $\ker\psi$ first, and then
prove that $[\ker\psi, \stl(m, n, \ml A)]=0$.

Denote by $P$ (resp. $Q$) the $K$-submodule of $\stl(m ,n, \ml A)$
generated by $v_{ij}(a)$ with $i<j$ (resp.  $i>j$). Clearly, we
have: $P=P_{0}\oplus P_{1}$ and $Q=Q_{0}\oplus Q_{1}$; the
restrictions of $\psi$ to $P$ and $Q$ are injective; the images of
$v_{ij}(a)$ from $P$ (resp. Q) under $\psi$ are strictly
uppertriangular (resp. lowertriangular) matrices in $\fsl(m, n,
\ml A)$.

Let $H_{ij}(a, b):=[v_{ij}(a), v_{ji}(b)]$ for $1\le i\ne j\le
m+n, a, b\in\ml A$,  $H$ the submodule of $\stl(m, n, \ml A)$
generated by $H_{ij}(a, b), i\ne j, a, b\in\ml A$. Then the
following Lemma is evident.

\begin{lemm}
Every element $X\in \stl(m, n, \ml A)$ with $m+n\ge 3$ can be uniquely written in the form
$$X=p+h+q, \quad\hbox{where}\ p\in P, \ h\in H, \ q\in Q.\eqno(3.1)$$
\end{lemm}

\begin{lemm}
For $m+n\ge3$, $\ker\psi\subseteq H$.
\end{lemm}
\begin{proof} Let $X=p+h+q\in\ker\psi$, where $ p\in P, \ h\in H, \ q\in Q$.
Then $0=\psi(X)=\psi(p)+\psi(h)+\psi(q)$. By Lemma 3.2, we have
$\psi(p)=\psi(q)=0$, so $p=q=0$, that is, $X=h\in H$.
\end{proof}

\begin{lemm}
For $m+n\ge3$, $[\ker\psi, \stl(m, n, \ml A)]=0$.
\end{lemm}
\begin{proof} By Lemma 3.3, any $t\in\ker\psi$ is expressible by
$H_{kl}(b, c)$'s. By definition, $$[v_{ij}(a), H_{kl}(b, c)]\in
P+Q.$$ So we have $[v_{ij}(a), t]=p+q$, where $ p\in P, q\in Q$.
Thus $\psi(p+q)=\psi([v_{ij}(a), t])=[\psi(v_{ij}(a)), \psi(t)]=0$
since $\psi(t)=0$. By injectivity of the restriction of $\psi$ to
$P+Q$, we get $p+q=0$. So $[\ker\psi, \stl(m, n, \ml A)]=0$.
\end{proof}

Therefore, we complete the proof of Theorem 3.1.

We proceed to show that the central extension $(\stl(m, n, \ml A),
\psi)$ is universal.
\begin{theo}
Let $(W, \phi)$ denote a central extension of Leibniz superalgebra
$\fsl(m, n, \ml A)$, and $m+n\ge 3$ with $\ch K\ne 2$ if $m+n=4$,
$\ch K\ne 3$ if $m+n=3$. Then there exists a unique homomorphism
$\rho: \stl(m, n, \ml A)\to W$ such that $\phi\circ\rho=\psi$.
\end{theo}
\begin{proof} \ For $m+n=3$ and $\ch K\ne 3$: \ one can use the same method in
the proof of Theorem 5.18 in \cite{AG} to prove $\stl(m, n, \ml A)$
is centrally closed. The differences here are that we need set
$\mathcal M=v_{13}(\ml A)\oplus v_{31}(\ml A)\oplus v_{23}(\ml
A)\oplus v_{32}(\ml A)$, $D=\ad(v_{12}(1)+v_{21}(1))$, and
$D|_{\mathcal M}$ is diagonalizable with eigenvalues $\pm 1$.

For  $m+n\ge 4$: \ since $\phi: W\longrightarrow \fsl(m, n, \ml A)$
is surjective, for any generator $E_{ij}(a)\in\fsl(m, n, \ml A)$, we
choose $e_{ij}(a)\in\phi^{-1}(E_{ij}(a))$. Then the commutator
$[e_{ij}(a), e_{kl}(b)]$ does not depend on the choice of
representatives of $\phi^{-1}(E_{ij}(a))$ and
$\phi^{-1}(E_{kl}(b))$. Moreover, for any $j\ne k$ and $i\ne l$, and
$a, \;b\in \ml A$, $\phi([e_{ij}(a), e_{kl}(b)])=[\phi(e_{ij}(a)),
\phi( e_{kl}(b))]=0$. So $[e_{ij}(a), e_{kl}(b)]\in \ker\phi$.

For distinct $i, j, k$, let
$$[e_{ik}(a), e_{kj}(b)]=e_{ij}(ab)+C_{ij}^k(a, b),
\quad [e_{kj}(b), e_{ik}(a)]=-(-1)^{\tau_{ik}\tau_{kj}}e_{ij}(ab)+D_{ij}^k(a, b),
$$ where $C_{ij}^k(a, b), \ D_{ij}^k(a, b)\in \ker\phi$.
Take $l\not\in\{i, j, k\}$, then
\begin{eqnarray*}
[e_{ik}(a), e_{kj}(bc)]&=&[e_{ik}(a), [e_{kl}(b), e_{lj}(c)]-C_{kj}^l(b, c)]=[e_{ik}(a), [e_{kl}(b), e_{lj}(c)]]\\
&=&[[e_{ik}(a), e_{kl}(b)], e_{lj}(c)]+(-1)^{\tau_{ik}\tau_{kl}}[e_{kl}(b), [e_{ik}(a), e_{lj}(c)]]\\
&=&[e_{il}(ab), e_{lj}(c)],
\end{eqnarray*}
i.e., $$[e_{ik}(a), e_{kj}(bc)]=[e_{il}(ab),
e_{lj}(c)].\eqno(3.2)$$

In particular,
$$[e_{ik}(a), e_{kj}(c)]=[e_{il}(a), e_{lj}(c)].\eqno(3.3)$$
It follows that $C_{ij}^k(a, c)=C_{ij}^l(a, c)$, which shows that
$C_{ij}^k$ is independent of the choice of $k$. Setting $C_{ij}(a,
b)=C_{ij}^k(a, b)$, we have
$$[e_{ik}(a), e_{kj}(b)]=e_{ij}(ab)+C_{ij}(a, b),\eqno(3.4)$$ where $C_{ij}(a, b)\in \ker\phi$.
Taking $a=1$, we have
$$[e_{ik}(1), e_{kj}(b)]=e_{ij}(b)+C_{ij}(1, b).\eqno(3.5)$$

Set $w_{ij}(a)=e_{ij}(a)+\ker\phi=[e_{ik}(1),
e_{kj}(a)]+\ker\phi$. We shall show that for all $i\ne j$, $a\in
\ml A$, $w_{ij}(a)$'s satisfy relations (1)---(4) in the
definition of Steinberg Leibniz superalgebra, i.e.,
\begin{gather*}
w_{ij}(k_1a+k_2b)=k_1w_{ij}(a)+k_2w_{ij}(b),\  \hbox{ for } \ a, \; b\in\ml A, \; k_1, \; k_2\in K;\tag{5}\\
[w_{ij}(a), w_{kl}(b)]=0,\  \hbox{ if }\ i\ne l, \  j\ne k;\tag{6}\\
[w_{ij}(a), w_{kl}(b)]=w_{il}(ab),\  \hbox{ if }\ i\ne l, \  j= k;\tag{7}\\
[w_{ij}(a), w_{kl}(b)]=-(-1)^{\tau_{ij}\tau_{kl}}w_{kj}(ba),\ \hbox{
if } \ i= l, \  j\ne k.\tag{8}
\end{gather*}

If $m+n\ge 5$, the proof of (5)---(8) is essentially the same as
that in \cite{KL} (or \cite{Gao2}, \cite{MP}).

If $m+n=4$, we only prove (6) for the case when $\{i, j, l,
k\}=\{1, 2, 3, 4\}$,
\begin{eqnarray*}
&&[w_{ij}(a), w_{kl}(bc)]=[w_{ij}(a), [w_{ki}(b), w_{il}(c)]]\\
&=&[[w_{ij}(a), w_{ki}(b)], w_{il}(c)]=-(-1)^{\tau_{ij}\tau_{ik}}[w_{kj}(ba), w_{il}(c)],\\
&&[w_{ij}(a), w_{kl}(bc)]=[w_{ij}(a), [w_{kj}(b), w_{jl}(c)]]\\
&=&(-1)^{\tau_{ij}\tau_{kj}}[w_{kj}(b), [w_{ij}(a),
w_{jl}(c)]]=(-1)^{\tau_{ij}\tau_{kj}}[w_{kj}(b), w_{il}(ac)].
\end{eqnarray*}
So we have
$$[w_{ij}(a), w_{kl}(bc)]=-(-1)^{\tau_{ij}\tau_{ik}}[w_{kj}(ba), w_{il}(c)]
=(-1)^{\tau_{ij}\tau_{kj}}[w_{kj}(b), w_{il}(ac)].\eqno(3.6)$$
Taking $a=b=1$ in (3.6), we have
$$[w_{ij}(1), w_{kl}(c)]=-[w_{kj}(1), w_{il}(c)]
=[w_{kj}(1), w_{il}(c)]=0,\eqno(3.7)$$ unless  $\tau_{ij}=1$
(noting
$(-1)^{\tau_{ij}\tau_{ik}}+(-1)^{\tau_{ij}\tau_{kj}}=(-1)^{\tau_{ij}\tau_{ik}}(1+(-1)^{\tau_{ij}})$\,).

Even in the case when $\tau_{ij}=1$, one can choose a
$k\not\in\{\,i,\;j\,\}$ such that $\tau_{kj}\ne1$ due to (2.1),
then use (3.7) again to obtain
$$[w_{ij}(1), w_{kl}(c)]=\pm[w_{kj}(1), w_{il}(c)]=0,\eqno(3.8)$$ for all  $\{i, j, l, k\}=\{1, 2, 3, 4\}$.

Taking $b=1$ in (3.6), together with (3.8), we obtain (6).

The others are also essentially the same as those in the case of
$m+n\ge 5$.

Now define a homomorphism $\rho: \stl(m, n, \ml A)\to W$  by
$\rho(v_{ij}(a))=w_{ij}(a)$. From the above, we see that this
mapping is well-defined, and $\phi\circ\rho=\psi$. The uniqueness of
the mapping $\rho$ follows from Lemmas 2.2, 2.3.
\end{proof}

\begin{remark} \  In (3.7), we used the restriction $\ch K\ne 2$.
From the proof of Theorem 5.18 in \cite{AG}, we must restrict $\ch
K\ne 3$ if $m+n=3$ (also see \cite{Gao1}).
\end{remark}

\section{Kernel of universal central extension $(\stl(m, n, \ml A), \psi)$}

In treating with the Leibniz superalgebras' case, in order to
calculate the kernel of the universal central extension of $\fsl(m,
n, \ml A)$, we have to introduce a so-called modified Hochschild
homology. By calculation, the modified Hochschild boundary must be
the $K$-linear map $d_n:\ml A^{\ot(n+1)}\longrightarrow \ml A^{\ot
n}$ defined by the formula (with the last summand different from the
usual definition of the Hochschild boundary, see \cite{Lo2}):
$$d_n(a_0\ot a_1\ot \cdots\ot a_n)
=\sum_{i=0}^{n-1}(-1)^ia_0\ot a_1\ot \cdots\ot
a_{i}a_{i+1}\ot\cdots\ot a_n-a_1\ot \cdots\ot a_{n-1}\ot a_na_0,$$
which still can be checked to satisfy $d_n\circ d_{n+1}=0$. Thus
one can consider the $n$-th homology group $HH_n(\ml A)=\ker
d_n/\im d_{n+1}$.

\begin{remark} \ The last summand in the definition (cf. \cite{Lo2}) of the
usual Hochschild boundary $d_n$ is $(-1)^na_na_0\ot a_1\ot \cdots\ot
a_{n-1}$.
\end{remark}

Let $M_{m+n}(\ml A)$ be the $K$-algebra of $(m+n)\times (m+n)$-matrices
written in the block form
$$\left(\begin{array}{c c} A & B \\  C & D \end{array}\right),$$ where
$A, B, C, D$ are matrices
in size $m\times m, m\times n, n\times m, n\times n$ respectively
with coefficients in $\ml A$, and $m, \ n\ge0$, $m+n\ge2$.

\begin{theo}
The kernel $T(m, n)$ of universal central extension $(\stl(m, n, \ml
A)$, $\psi)$ of $\fsl(m, n, \ml A)$ in {\bf SLeib} under the
assumption $m+n\ge 3$ $($with $\ch K\ne 2$ if $m+n=4$, $\ch K\ne 3$
if $m+n=3$$)$ is isomorphic to $HH_1(\ml A)$.
\end{theo}
\begin{proof} \  The proofs are analogous to those in \cite{KL}
(also see \cite{Gao2}, \cite{MP}), except for some properties of
$H_{ij}(a, b)$ and the homomorphism $\theta: \stl(m, n, \ml A)
\longrightarrow \ml A\ot\ml A/\im d_2 $. In our case, $H_{ij}$'s
satisfy the following properties:

\smallskip
$1)$ $H_{ij}(ab, c)=H_{ik}(a, bc)+(-1)^{\tau_{ik}}H_{kj}(b, ca)$;

\smallskip
$2)$ $H_{ij}(1, a)=-(-1)^{\tau_{ij}}H_{ji}(1, a)$, for any
distinct $i, j$ (we also set $h(a, b)=H_{1i}(a, b)-H_{1i}(1, ba)$,
which is independent of $i(\ne 1$)).

\smallskip
$\theta$ is defined by $\theta([x, y])=\sum_{i, j}\psi(x)_{ij}\ot
\psi(y)_{ji}=\str_2(\psi(x)\ot\psi(y))$, where $\str_2:M_{m+n}(\ml
A)\ot M_{m+n}(\ml A)\longrightarrow \ml A\ot\ml A$ is given by
\begin{eqnarray*}
&&\str_2(\left(\begin{array}{c c} A_1 & B_1 \\  C_1 & D_1 \end{array}\right)\ot \left(\begin{array}{c c} A_2 & B_2 \\  C_2 & D_2 \end{array}\right))\\
&=&\tr_2(\left(\begin{array}{c c} A_1 & 0 \\  0 &0
\end{array}\right)\ot \left(\begin{array}{c c} A_2 & 0 \\ 0 & 0
\end{array}\right)+
\left(\begin{array}{c c} 0 & B_1 \\ 0 & 0 \end{array}\right)\ot \left(\begin{array}{c c} 0 & 0 \\  C_2 & 0 \end{array}\right))\\
&&\ -\,\tr_2(\left(\begin{array}{c c} 0 & 0 \\  C_1 & 0
\end{array}\right)\ot \left(\begin{array}{c c} 0 & B_2 \\  0 & 0
\end{array}\right)+ \left(\begin{array}{c c} 0 & 0 \\  0 & D_1
\end{array}\right)\ot \left(\begin{array}{c c} 0 & 0 \\  0 & D_2
\end{array}\right)),
\end{eqnarray*}
and $\tr_2$ was defined in \cite{KL} by $\tr_2(P\ot Q)=\sum_{1\le i,
j\le l}p_{ij}\ot q_{ji}$, for $P, \, Q\in M_l(\ml A)$.
\end{proof}

\begin{coro} $($\cite{LH3}$)$ \
When $\ml A$ is commutative, $HH_1(\ml A)\cong\Omega_{\ml A}^1$,
i.e.,
$$0\longrightarrow \Omega_{\ml A}^1\longrightarrow \stl(m, n, \ml A)\longrightarrow \fsl(m , n, \ml A)
\longrightarrow 0$$ is the universal central extension of $\fsl(m,
n, \ml A)$ in {\bf SLeib} under the assumption: $m+n\ge 3$ with $\ch
K\ne 2$ if $m+n=4$, or $\ch K\ne 3$ if $m+n=3$.
\end{coro}

\begin{remark} \ If take $n=0$ in $\fsl(m, n, \ml A)$, we recover the same
results in \cite{LP}, \cite{Gao2} for Leibniz algebras. The stable
case (characteristic 0) is due to Cuvier-Loday theorem (see
\cite{Lo2}).
\end{remark}

\section{Universal central extension of $\fst(m, n, \ml A)$ in {\bf SLeib}}

In what follows, it is necessary for us to point out an important
relationship between our modified Hochschild homology (only
applicable to the Leibniz superalgebras' case) and the cyclic
homology.

Consider the complex of the $K$-modules $C_*(\ml A)$, where $C_0(\ml
A)=\ml A$ and for $n\ge1$ the module $C_n(\ml A)$ is the factor
module of the $K$-module $\ml A^{\ot(n+1)}$ by the $K$-submodule
generated by elements $a_0\ot\cdots\ot a_n-(-1)^na_1\ot\cdots\ot
a_n\ot a_0, a_i\in \ml A, i=0, \cdots, n$. The cyclic boundary
induced by the modified Hochschild boundary (also denoted $d_n$) is
exactly the one induced by the usual Hochschild boundary.  Thus we
consider the essentially same $n$-th cyclic homology group $HC_n(\ml
A)=\ker d_n/\im d_{n+1}$ as usual in \cite{Lo2}. Moreover, there is
a natural projection $p: HH_*(\ml A)\longrightarrow HC_*(A)$ induced
by the projection $\ml A^{\ot(n+1)}\longrightarrow C_n(\ml A)$. Note
the well-known exact sequence
$$0\longrightarrow HC_1(\ml A)\longrightarrow C_1(\ml A)/\im d_2\stackrel{d_1}{\to} \ml A,$$
as well as the Connes operator $\ml B$ (see \cite{Lo2}), which is a
$K$-linear map $\ml A^{\ot(n+1)}\longrightarrow \ml A^{\ot(n+2)}$ by
\begin{equation*}
\begin{split}
 \ml B(a_0\ot\cdots\ot a_n)&=\sum_{i=0}^n(-1)^{ni}(1\ot
a_i\ot\cdots\ot a_n\ot a_0\ot\cdots a_{i-1})\\
&\quad+(-1)^{n(i+1)}(a_i\ot\cdots\ot a_n\ot a_0\ot\cdots a_{i-1}\ot
1),
\end{split}
\end{equation*}
such that $\ml Bd_n+d_{n+1}\ml B=0$. Using these objects,
one can prove the following theorem.
\begin{theo} \cite{MP}
 The universal central extension of $\fsl(m, n, \ml
A)$ in category {\bf SLie} is $\fst(m, n, \ml A)$ with kernel
$HC_1(\ml A)$ under the assumption: $m+n\ge5$.
\end{theo}

\begin{remark} \ Using the same methods in the proofs of
Sections 3 \& 4, we can show that Theorem 5.1 still holds under a
weak assumption below.
\end{remark}
\begin{theo}
The universal central extension of $\fsl(m, n, \ml A)$ in category
{\bf SLie} is $\fst(m, n, \ml A)$ with kernel $HC_1(\ml A)$ under
the assumption: $m+n\ge 3$ with $\ch K\ne 2$ if $m+n=4$, $\ch K\ne
3$ if $m+n=3$.
\end{theo}

Now we give the main theorem of this section.
\begin{theo} Steinberg Leibinz superalgebra $(\stl(m, n, \ml A), \pi)$ is
the universal central extension of $\fst(m, n, \ml A)$ in category
{\bf SLeib} with the kernel isomorphic to $\im\ml B$ under the
assumption: $m+n\ge 3$ with $\ch K\ne 2$ if $m+n=4$, or $\ch K\ne 3$
if $m+n=3$.
\end{theo}

\begin{proof} \ The proving idea is similar to that in \cite{LP}.

Let $\mg=\stl(m, n, \ml A)$, then $\mg_{SLie}=\fst(m, n, \ml A)$.
From \cite{MP} or Theorem 5.3, we see that $\mg_{SLie}$ is the
universal central extension of $\fsl(m, n, \ml A)$ in {\bf SLie} for
$m+n\ge3$. By Theorem 3.5, Theorem 4.2 and Theorem 5.3, we have the
following exact diagram. Moreover,  it is clear that this diagram is
commutative.
$$\begin{array}{ c c c c c c c }
           & HC_0(\ml A)&&&&&\\
           &\quad\downarrow\ml B&&&&& \\
           0\longrightarrow &HH_1(\ml A)&\longrightarrow&\stl(m, n, \ml A)&\stackrel{\psi}{\longrightarrow}&\fsl(m, n, \ml A)&\longrightarrow0 \\
        &\quad\downarrow p&&\quad\downarrow \pi&&\parallel& \\
           0\longrightarrow &HC_1(\ml A)&\longrightarrow&\fst(m, n, \ml A)&\stackrel{\vp}{\longrightarrow}&\fsl(m, n, \ml A)&\longrightarrow0\\
        &\downarrow &&\downarrow&&&\\
        &0&&0&&&
\end{array}$$
Then we have
$$\ker(\pi)\cong \ker p\cong \im \ml B.$$
Hence $0\longrightarrow \im\ml B\longrightarrow \stl(m, n, \ml
A)\longrightarrow \fst(m, n, \ml A)\longrightarrow 0$ is a central
extension in {\bf SLeib}. Moreover, for $m+n\ge3$, this is a
universal extension by Theorem 5.3 and Theorem 4.2. Therefore,
$$HL_2(\fst(m, n, \ml A))\cong\im\ml B.  $$
\end{proof}

\begin{remark} \ Taking $n=0$ in Theorem 5.4, we obtain that (4.6) in \cite{LP}
also holds for all $m\ge3$ with $\ch K\ne 2$ if $m=4$, $\ch K\ne 3$
if $m=3$ (see \cite{Gao2}).
\end{remark}

\bigskip
\section*{Acknowledgments}
\medskip

N. Hu thanks N. Bergeron, Y. Gao for financial support during a
one-year visit in the Department of Mathematics and Statistics of
York University. The part of this work was done during D. Liu's
Ph.~D study in the Department of Mathematics at University of
Bielefeld under the support of the `Asia-Euro Link Programme'. He
expresses his deep thanks to C.M. Ringel for the continuous
encouragement. Authors give their special thanks to Gao for the
crucial comments.

\bigskip
\bibliographystyle{amsalpha}

\end{document}